\documentstyle[12pt]{article}
\newtheorem{theo}{Theorem}[section]

\newtheorem{coro}[theo]{Corollary}
\newtheorem{rema}[theo]{Remark}
\newtheorem{prop}[theo]{Proposition}
\def\al{\alpha}
\def\be{\beta}

\def\de{\delta}

\def\th{\theta}

\def\rh{\rho}
\def\si{\sigma}

\def\ph{\phi}
\def\vp{\varphi}

\def\De{\Delta}
\def\La{\Lambda}

\def\Box{{\vbox{\hrule\hbox{\vrule\phantom{s}\vrule}\hrule}\,}}
\def\beq{\begin{equation}}
\def\eeq{\end{equation}}
\def\bea{\begin{eqnarray}}
\def\eea{\end{eqnarray}}
\def\beas{\begin{eqnarray*}}
\def\eeas{\end{eqnarray*}}
\def\nn{\nonumber}

\font\twlmsbm=msbm10 at 12 true pt

\font\egtmsbm=msbm8
\font\sevmsbm=msbm7
\newfam\msbmfam
\textfont\msbmfam=\twlmsbm
\scriptfont\msbmfam=\egtmsbm
\scriptscriptfont\msbmfam=\sevmsbm
\def\Bbb#1{{\fam\msbmfam\relax#1}} 
\def\Zah{\Bbb Z}
\def\Real{\Bbb R}

\def\C{\Bbb C}
\font\twleufm=eufm10 at 12 true pt

\font\egteufm=eufm8
\font\seveufm=eufm7
\newfam\eufmfam
\textfont\eufmfam=\twleufm
\scriptfont\eufmfam=\egteufm
\scriptscriptfont\eufmfam=\seveufm
\def\Eu#1{{\fam\eufmfam\relax#1}} 
\def\su{{\Eu s \Eu u}}
\def\U{{\cal U}_q(\su(1,1))}
\begin{document}
\begin{center}
{\large \bf BILINEAR GENERATING FUNCTIONS}\\[3mm]
{\large \bf FOR ORTHOGONAL POLYNOMIALS}\\[2cm]
{\sc H.T.\ Koelink\footnote{Supported by the Netherlands Organization
for Scientific Research (NWO) under project number 610.06.100.} and 
J.~Van der Jeugt}\footnote{Research Associate of the Fund
for Scientific Research -- Flanders (Belgium)} \\[1cm]
Report 97-03, Universiteit van Amsterdam\\[5mm]
April 24, 1997 \\[2mm]
\end{center}

\begin{abstract}
Using realisations of the positive discrete series representations of
the Lie algebra $\su(1,1)$ in terms of Meixner-Pollaczek
polynomials, the action of $\su(1,1)$ on Poisson kernels of
these polynomials is considered. In the tensor product of two such
representations, two sets of eigenfunctions of a certain operator can
be considered and they are shown to be related through
continuous Hahn polynomials. As a result, a bilinear generating
function for continuous Hahn polynomials is
obtained involving the Poisson kernel of
Meixner-Pollaczek polynomials. For the positive discrete series
representations of the quantised universal enveloping algebra $\U$ a
similar analysis is performed and leads to a bilinear generating
function for Askey-Wilson polynomials 
involving the Poisson kernel of Al-Salam and Chihara polynomials.
\end{abstract}

\vskip 1cm
{\small
\noindent 1991 Mathematics Subject Classification~: 33C80, 33D80,
33C25, 33D25, 17B20, 17B37.\\[1mm]
\noindent Keywords~: orthogonal polynomials, bilinear generating
functions, Poisson kernels, Lie algebra, quantum algebra.
}
\vskip 1cm
%
\section{Introduction}  \label{sec-intro}

Representation theory of Lie algebras and quantum algebras (or
quantised universal enveloping algebras~\cite{CharP}) is closely
related to special functions of (basic) hypergeometric
type~\cite[Ch.~14]{VileK}, \cite[Ch.~13]{CharP}. 
In this paper we are dealing with positive
discrete series representations of the Lie algebra $\su(1,1)$ and of its
quantum analogue $\U$. Using realisations of these representations in
terms of orthogonal polynomials $p_n(x)$, the action of the Lie or quantum
algebra on Poisson kernels of these polynomials is considered. In the
realisation of the tensor product of two representations, one can
either use the product $p_{n_1}(x_1) p_{n_2}(x_2)$ as basis vectors, or
a new basis in which different orthogonal polynomials $q_n(x)$
appear~\cite{KV}. Here it is shown how bilinear generating functions
for the $q_n(x)$ appear naturally in this framework. 

For the Lie algebra $\su(1,1)$ it was shown in~\cite{V,KV}, using ideas
of Granovskii and Zhedanov~\cite{GZ}, that
Laguerre, Meixner-Pollaczek and Meixner polynomials appear as overlap
coefficients in the positive discrete series representations. As a
consequence, basis vectors of these representations have a realisation
in terms of these polynomials. In the realisation with Laguerre
polynomials, the tensor product of two representations has a so-called
uncoupled basis with basis vectors the product of two Laguerre
polynomials, or a coupled basis with basis vectors the product of a
Jacobi and a Laguerre polynomial. The Clebsch-Gordan coefficients of
$\su(1,1)$ relate one basis to another, and as a result a new
convolution formula for Laguerre polynomials involving a Jacobi
polynomial follows~\cite{V}. This result was extended yielding a new
convolution formula for Meixner-Pollaczek (resp.\ Meixner) polynomials
involving a continuous Hahn (resp.\ Hahn) polynomial~\cite{KV}. 

In the present
paper we use the action of the Lie algebra on Poisson kernels. Let us
briefly describe the situation and the technique when the
realisation is in terms of Meixner-Pollaczek polynomials. Then a
positive discrete series representation is identified with
$L^2(\Real,d\mu(x))$, $d\mu(x)$ being the orthogonality measure of
Meixner-Pollaczek polynomials, and its basis vectors are realised in terms
of Meixner-Pollaczek polynomials. 
The Lie algebra $\su(1,1)$ has a realisation in
terms of operators acting in $L^2(\Real,d\mu(x))$. The
Poisson kernel can be seen as an expansion of some function in terms of
Meixner-Pollaczek polynomials, thus $\su(1,1)$ has a natural action on it. In
this case the Poisson kernel is proportional to a ${}_2F_1$-series.
The Poisson kernel is shown to be an
eigenfunction of some element $X_t$ (in this realisation) of $\su(1,1)$.
This element $X_t$ is closely related to a recurrence operator that led
to the interpretation of Meixner-Pollaczek polynomials as overlap coefficients
in the positive discrete series representations. Then we go on
considering eigenfunctions of $X_t$ in the tensor product of two
representations. There are again so-called uncoupled eigenfunctions,
being the product of two Poisson kernels, and coupled eigenfunctions.
The Clebsch-Gordan coefficients relating these two are the polynomials
that appear in the convolution formula, i.e.\ continuous Hahn polynomials.
As a consequence of this relation, one obtains a bilinear generating
function for continuous Hahn polynomials~(Theorem~\ref{theochahn}),
involving the Poisson kernel of
Meixner-Pollaczek polynomials. This appears to be a new formula, which
also has an interpretation as a multiplication formula for
hypergeometric series~(Theorem~\ref{theoMP}).
Theorem~\ref{theochahn}
generalises Bateman's 1904 bilinear generating function for the Jacobi
polynomials. The same technique can be applied to the realisation
in terms of Meixner polynomials, for which we give only the final result.

In section~3 we show how this approach can be applied to the
quantised universal enveloping algebra $\U$. The basis vectors of the
positive discrete series representations of $\U$ can be realised in
terms of Al-Salam and Chihara polynomials~\cite{KV}. The general
convolution identity for Al-Salam and Chihara polynomials,
see~\cite[Theorem~4.5]{KV}, involves an Askey-Wilson polynomial. For
Al-Salam and Chihara polynomials the Poisson kernel needed here was
given recently in a paper by Askey, Rahman and Suslov~\cite{ARS}; another
method by which this Poisson kernel can be obtained is given
in~\cite{JJ}. This Poisson kernel is expressed by means of a
very-well-poised basic hypergeometric series ${}_8\vp_7$, often denoted
by ${}_8W_7$, where we follow the notation of Gasper and Rahman~\cite{GR}
for basic hypergeometric series. 
In the present case, our technique leads to a bilinear
generating function for Askey-Wilson polynomials; the only
other functions appearing in this formula being ${}_8W_7$
series~(Theorem~\ref{theoAW}). Limiting cases yield a
bilinear generating function for continuous dual $q$-Hahn polynomials
and for Al-Salam and Chihara polynomials.
The bilinear generating function for Askey-Wilson polynomials
(Theorem~\ref{theoAW}) is a $q$-analogue of Bateman's bilinear
generating function for Jacobi polynomials.

\newfont{\eultit}{eufm10 scaled\magstep2} 
\section{Polynomials related to $\mbox{\eultit su}(1,1)$}
\setcounter{equation}{0}
\def\Zp{\Zah_+}
\def\Hi{\ell^2(\Zah_+)}

The Lie algebra $\su(1,1)$ is generated by $H,B,C$ subject to
$$
[H,B] = 2B, \qquad [H,C] = -2C, \qquad [B,C]=H.
$$
There is a $\ast$-structure by $H^\ast=H$ and $B^\ast=-C$.
The positive discrete series representations $\pi_k$ of
$\su(1,1)$ are unitary representations labelled by $k>0$.
The representation space is $\Hi$ equipped
with an orthonormal basis $\{ e^k_n\}_{\{ n\in\Zp\} }$, and the action
is given by
\beq
\begin{array}{ll}
\pi_k(H)\, e^k_n &= 2(k+n) \, e^k_n, \\[1mm]
\pi_k(B)\, e^k_n &= \sqrt{(n+1)(2k+n)}\, e^k_{n+1},\\[1mm]
\pi_k(C)\, e^k_n &= -\sqrt{n(2k+n-1)}\, e^k_{n-1}.
\end{array}
\label{actionsu11}
\eeq
The tensor product of two positive discrete series
representations decomposes as
\beq
\pi_{k_1}\otimes \pi_{k_2} = \bigoplus_{j=0}^\infty
\pi_{k_1+k_2+j},
\label{tensorproddec}
\eeq
and the corresponding intertwining operator can be expressed
by means of the Clebsch-Gordan coefficients
\beq
e^{(k_1k_2)k}_n = \sum_{n_1,n_2} C^{k_1,k_2,k}_{n_1,n_2,n}\,
e^{k_1}_{n_1}\otimes e^{k_2}_{n_2}.
\label{defCGC}
\eeq
The Clebsch-Gordan coefficients are non-zero only if $n_1+n_2=n+j$,
$k=k_1+k_2+j$ for $j,n_1,n_2,n\in\Zp$, and normalised
by $\langle e^k_0, e^{k_1}_0\otimes e^{k_2}_j\rangle >0$.
For the above results see Vilenkin and Klimyk~\cite[\S 8.7]{VileK},
or~\cite{KV}. 

Up to a scalar multiple, the most general self-adjoint element in
$\su(1,1)$ is of the form
\beq
X_\al=B-C-\al H,\qquad \al\in\Real,
\label{X}
\eeq
thus $X_\al^*=X_\al$. The action of $X_\al$ on a basis vector $e^k_n$ in the
representation $\pi_k$ is then given by
\bea
&&\pi_k(X_\al) e^k_n = a_{n} e^k_{n+1} + b_n e^k_n + a_{n-1} e^k_{n-1}
\label{actionX} \\
&&\hbox{where } a_n=\sqrt{(n+1)(2k+n)} \hbox{ and } b_n=-2\al(k+n). \nn
\eea
Thus $\pi_k(X_\al)$ acts as a three-term recurrence operator.
Since the coefficients $a_n$ and $b_n$ in
(\ref{actionX}) satisfy $a_n>0$ and $b_n\in\Real$, there exist
orthonormal polynomials associated with $\pi_k(X_\al)$ defined
by~\cite{Berez} 
\bea
&&x p_n(x) = a_{n} p_{n+1}(x) + b_n p_n(x) + a_{n-1} p_{n-1}(x)
\label{recurrence} \\
&&\hbox{with } p_{-1}(x)=0, \quad p_0(x)=1. \nn
\eea
In~\cite{KV} it was shown that these polynomials are Meixner-Pollaczek
polynomials if $|\al|<1$, Meixner polynomials if $|\al|>1$, and
Laguerre polynomials if $|\al|=1$ (this last case was also considered
in~\cite{V}). A consequence of this is that the basis vectors of the
positive discrete series representations can be realised as these
polynomials, and the elements of $\su(1,1)$ are realised as operators in
the corresponding $L^2$ space such that $X_\al$ is realised as a 
multiplication operator. 

Here the realisation in
terms of Meixner-Pollaczek polynomials is considered. It is shown how the
$\su(1,1)$ operators act on Poisson kernels of these polynomials, and
that they are eigenfunctions of some operator $X_t$. Then two
different eigenfunctions of this operator in the tensor product are
introduced, and their relationship is shown to lead to a bilinear
generating function for continuous Hahn polynomials. The case of Laguerre
polynomials (and related Jacobi polynomials) is obtained as a limit.
The case of Meixner polynomials
(and related Hahn polynomials) is completely analogous, and we just
mention the final result at the end of this paragraph.

The Meixner-Polla\-czek polynomials are defined by
\beq
P_n^{(k)}(x;\phi) = {{(2k)_n}\over{n!}} e^{in\phi}
\, {}_2F_1\left( {{-n,k+ix}\atop{2k}};1-e^{-2i\phi}\right),
\label{defMP}
\eeq
with the usual notation for Pochhammer symbols and hypergeometric 
series~\cite{GR}.
For $k>0$ and $0<\phi<\pi$ these are orthogonal polynomials
with respect to a positive measure on $\Real$, see \cite{KoekS}, 
\cite[App.]{Szego}.
The orthonormal Meixner-Pollaczek polynomials
$$
p_n^{(k)}(x;\phi) = \sqrt{ {{n!}\over{\Gamma(n+2k)}}}
P_n^{(k)}(x;\phi)
$$
satisfy the orthogonality relation
$$
\int_\Real p_n^{(k)}(x;\phi)p_m^{(k)}(x;\phi)\,
{{ (2\sin\phi)^{2k}}\over{2\pi}}
e^{(2\phi-\pi)x} \bigl\vert \Gamma(k+ix)\bigr\vert^2 dx =\de_{n,m}.
$$
We shall denote this absolutely continuous measure by $d\mu_{k,\ph}(x)$.
There exists~\cite[Proposition~3.1]{KV} a unitary operator 
$\La : \Hi \rightarrow L^2(\Real,
d\mu_{k,\phi}(x))$ mapping each $e^k_n$ onto $p^{(k)}_n(x;\phi)$, and
every element $X$ of $\su(1,1)$ has a realisation $\rh_{k,\phi}(X)$ in 
$L^2(\Real, d\mu_{k,\phi}(x))$.

\begin{prop} 
For $t\in\C$ with $|t|<1$ and $0<\phi<\pi$ denote
$$
r = -4t\sin^2\phi / (1-t)^2.
$$
Let $y\in\Real$ and define
\beq
v^k_t(x,y;\phi) = \sum_{n=0}^\infty p^{(k)}_n(x;\phi)
p^{(k)}_n(y;\phi) t^n. 
\label{vMP}
\eeq
Then\\
(i) $ v^k_t(x,y;\phi) \in L^2(\Real,d\mu_{k,\phi}(x))$.\\
(ii) Explicitly,
\beq
v^k_t(x,y;\phi) = {1\over{\Gamma(2k)}}
(1-te^{2i\phi})^{i(x+y)} 
(1-t)^{-2k-ix-iy}  
\, {}_2F_1\left( {{k+ix, k+iy}\atop{2k}} ; r \right).
\label{vexplicitMP}
\eeq
(iii) Let 
$$
X_t= -\cos\phi\, H + tB - t^{-1}C \in \su(1,1).
$$ 
Then $v^k_t(x,y;\phi)$ is an
eigenfunction of $\rh_{k,\phi}(X_t)$ for the eigenvalue $2y\sin\phi$.
\label{propMP}
\end{prop}

\noindent {\em Proof.}
The first statement is true when $\sum_n |p^{(k)}_n(y;\phi) t^n|^2 < 
\infty$, and this holds since $p^{(k)}_n(y;\phi)$ are polynomials of the 
first kind associated with a three-term recurrence 
operator~\cite[Ch. VII, (1.24)]{Berez}. 
Statement (i) also follows for $y\in\C$ from
the asymptotic behaviour of the Meixner-Pollaczek polynomials.
For $k=1/2$ see Szeg\H o \cite[App. \S 6]{Szego}, and the general
case follows similarly using Darboux's method \cite[\S 8.4]{Szego}
on \cite[App. (4.1)]{Szego}. 
The explicit formula (\ref{vexplicitMP})
is the Poisson kernel for the
Meixner-Pollaczek polynomials and follows from~\cite[2.5.2 (12)]{HTF1}.
The last statement follows from (\ref{actionsu11}) and the three-term
recurrence relation 
$$
2y\sin\phi\; p^{(k)}_n(y;\phi) = a_n p^{(k)}_{n+1}(y;\phi) -
2(n+k)\cos\phi\; p^{(k)}_n(y;\phi) + a_{n-1} p^{(k)}_{n-1}(y;\phi),
$$
where $a_n=\sqrt{(n+1)(n+2k)}$. \qquad\Box

Note that the nonsymmetric Poisson kernel $\sum_{n=0}^\infty
p^{(k)}_n(x;\phi) p^{(k)}_n(y;\psi) t^n$, i.e.\ (\ref{vMP}) with
$\phi\ne\psi$, can also be determined 
explicitly using~\cite[2.5.2 (12)]{HTF1}, but it does not lead to a
more general result.

In the tensor product $\pi_{k_1}\otimes \pi_{k_2}$,
there is a unitary operator $\Upsilon : \Hi\otimes\Hi \rightarrow L^2(
\Real^2, d\mu_{k_1,\phi_1}(x_1) d\mu_{k_2,\phi_2}(x_2))$ mapping
$e^{k_1}_{n_1} \otimes e^{k_2}_{n_2}$ onto $p^{(k_1)}_{n_1}(x_1;\phi_1)
p^{(k_2)}_{n_2}(x_2;\phi_2)$. We shall only consider the case $\phi_1 =
\phi_2 = \phi$, because then we can simplify the basis
functions in the tensor product decomposition. The realisation of
$\Delta(X)=X\otimes 1 + 1 \otimes X$ with $X\in \su(1,1)$ in $L^2(
\Real^2, d\mu_{k_1,\phi}(x_1) d\mu_{k_2,\phi}(x_2))$ is denoted by 
$\si_{k_1,k_2,\phi}(X)$.
For $k=k_1+k_2+j$, with $j\in\Zp$, the above-mentioned simplification
reads~\cite[(3.7)]{KV}
\bea
e^{(k_1k_2)k}_n(x_1,x_2;\phi) &=& \sum_{n_1,n_2} C^{k_1,k_2,k}_{n_1,n_2,n}\
 p^{(k_1)}_{n_1}(x_1;\phi)  p^{(k_2)}_{n_2}(x_2;\phi)  \nn\\
&=& p^{(k)}_n(x_1+x_2;\phi)\, S^{(k_1,k_2)}_j(x_1,x_2;\phi),
\label{ecoupledMP}
\eea
where 
\bea
S^{(k_1,k_2)}_j(x_1,x_2;\phi)&=& (-2\sin\phi)^j
\sqrt{ {{j!\, (2j+2k_1+2k_2-1)\Gamma(j+2k_1+2k_2-1)}\over
{\Gamma(2k_1+j) \Gamma(2k_2+j)}} } \nn\\
&&\times p_j(x_1;k_1, k_2-i(x_1+x_2),k_1,k_2+i(x_1+x_2)).
\label{SexplicitMP}
\eea
Herein, $p_j$ is a continuous Hahn polynomial introduced by Atakishiyev
and Suslov~\cite{AS}, see also~\cite{Askey,KoekS},
\beq
p_n(x;a,b,c,d)= i^n {{(a+c)_n(a+d)_n}\over{n!}}\,
{}_3F_2\left( {{-n,n+a+b+c+d-1,a+ix}\atop{a+c,\ a+d}};1\right).
\label{defchahn}
\eeq

Next, we have

\begin{prop}
Let $y,t$ and $\phi$ be as in Proposition~\ref{propMP}, and
\beq
w^{(k_1k_2)k}_t(x_1,x_2,y;\phi) = \sum_{n=0}^\infty p^{(k)}_n(y;\phi)\,
 e^{(k_1k_2)k}_n(x_1,x_2;\phi)\, t^n
\label{wMP}
\eeq
for $k=k_1+k_2+j$ with $j\in\Zp$.
Let again $X_t=-\cos\phi\, H + tB - t^{-1}C\in \su(1,1)$, then\\
(i) $w^{(k_1k_2)k}_t(x_1,x_2,y;\phi) \in L^2(\Real^2,
d\mu_{k_1,\phi}(x_1) d\mu_{k_2,\phi}(x_2))$ and moreover it is an
eigenfunction of $\si_{k_1,k_2,\phi} (X_t)$ for the eigenvalue
$2y\sin\phi$.\\ 
(ii) Explicitly, one has 
\beq
w^{(k_1k_2)k}_t(x_1,x_2,y;\phi) = S_j^{(k_1,k_2)}(x_1,x_2;\phi)\,
v^k_t(x_1+x_2,y;\phi). 
\label{wexplicitMP}
\eeq
(iii) For $y_1,y_2\in\Real$,
\beq
v^{k_1}_t(x_1,y_1;\phi)v^{k_2}_t(x_2,y_2;\phi) =
\sum_{j=0}^\infty t^j S_j^{(k_1,k_2)}(y_1,y_2;\phi)\,
w^{(k_1k_2)k}_t(x_1,x_2,y_1+y_2; \phi).
\label{CGwMP}
\eeq
\label{propwMP}
\end{prop}

\noindent
The proposition states that in $L^2(\Real^2,d\mu_{k_1,\phi}(x_1)
d\mu_{k_2,\phi}(x_2))$ the uncoupled eigenfunc\-tions of 
$\si_{k_1,k_2,\phi}(X_t)$ for
the eigenvalue $2(y_1+y_2)\sin\phi$ are $v^{k_1}_t(x_1,y_1;\phi)
v^{k_2}_t(x_2,y_2;\phi)$;
the coupled eigenfunctions for the same eigenvalue are
$w^{(k_1k_2)k}_t(x_1,x_2,y_1+y_2;\phi)$; and the Clebsch-Gordan coefficients
relating these two are $t^j S_j^{(k_1,k_2)}(y_1,y_2;\phi)$. \\[2mm]
{\em Proof.} Since the action of $\si_{k_1,k_2,\phi}(H)$, 
$\si_{k_1,k_2,\phi}(B)$
and $\si_{k_1,k_2,\phi}(C)$ on the basis functions 
$e^{(k_1k_2)k}_n(x_1,x_2;\phi) \in
L^2(\Real^2,d\mu_{k_1,\phi}(x_1) d\mu_{k_2,\phi}(x_2))$ is 
the standard one given
by (\ref{actionsu11}), one can use again the three-term recurrence
relation for 
the orthonormal Meixner-Pollaczek polynomials to see that 
$\si_{k_1,k_2,\phi}(X_t)
w^{(k_1k_2)k}_t(x_1,x_2,y;\phi) = 2y\sin\phi\; 
w^{(k_1k_2)k}_t(x_1,x_2,y;\phi)$. 
Using (\ref{ecoupledMP}) in (\ref{wMP}), and then definition
(\ref{vMP}), (\ref{wexplicitMP}) follows directly. Finally, from
(\ref{vMP}) 
\bea
&&v^{k_1}_t(x_1,y_1;\phi)v^{k_2}_t(x_2,y_2;\phi) = \nn\\
&&\sum_{n_1=0}^\infty p^{(k_1)}_{n_1}(x_1;\phi) p^{(k_1)}_{n_1}(y_1;\phi)\,
 t^{n_1} \sum_{n_2=0}^\infty p^{(k_2)}_{n_2}(x_2;\phi) 
p^{(k_2)}_{n_2}(y_2;\phi)\, t^{n_2}.
\label{vv1}
\eea
Herein, use
$$
p^{(k_1)}_{n_1}(x_1;\phi) p^{(k_2)}_{n_2}(x_2;\phi) = \sum_{k,n}
C^{k_1,k_2,k}_{n_1,n_2,n}\, e^{(k_1k_2)k}_n(x_1,x_2;\phi).
$$
This equation follows from (\ref{ecoupledMP}) and the orthogonality of
the Clebsch-Gordan coefficients; the finite sum is over all $k=k_1+k_2+j$
($j\in\Zp$) and $n$ such that $j+n=n_1+n_2$. Then (\ref{vv1}) reduces
to 
\beas
&&v^{k_1}_t(x_1,y_1;\phi)v^{k_2}_t(x_2,y_2;\phi) \nn\\
&&= \sum_{k,n} t^{n+j}
e^{(k_1k_2)k}_n(x_1, x_2;\phi) \sum_{n_1,n_2} C^{k_1,k_2,k}_{n_1,n_2,n}\,
p^{(k_1)}_{n_1}(y_1;\phi) p^{(k_2)}_{n_2}(y_2;\phi) \\
&&= \sum_{k,n} t^{n+j} e^{(k_1k_2)k}_n(x_1, x_2;\phi) 
p^{(k)}_n(y_1+y_2;\phi) S^{(k_1,k_2)}_j(y_1,y_2;\phi) \\
&&= \sum_{j=0}^\infty t^j S^{(k_1,k_2)}_j(y_1,y_2;\phi) 
w^{(k_1k_2)k}_t(x_1,x_2,y_1+y_2;\phi),
\eeas
proving the final statement. All series manipulations above are allowed
because of the absolute convergence of the power series in $t$
involved.  \qquad\Box

The basic type of identity in this paper follows from 
(\ref{wexplicitMP}) and (\ref{CGwMP})~: 
\bea
&&v^{k_1}_t(x_1,y_1;\phi)v^{k_2}_t(x_2,y_2;\phi) =
\label{SpoissonMP} \\
&&\sum_{j=0}^\infty t^j v^{k_1+k_2+j}_t(x_1+x_2,y_1+y_2;\phi)
S_j^{(k_1,k_2)}(x_1,x_2;\phi)\,S_j^{(k_1,k_2)}(y_1,y_2;\phi). \nn
\eea
Using (\ref{vexplicitMP}) and (\ref{SexplicitMP})
some simplifications take place, and (\ref{SpoissonMP}) becomes
($k=k_1+k_2+j$) 
\bea
&&\, {}_2F_1\left( {{k_1+ix_1, k_1+iy_1}\atop{2k_1}};r \right)
{\ }_2F_1\left( {{k_2+ix_2, k_2+iy_2}\atop{2k_2}}; r\right) = \nn\\
&&\sum_{j=0}^\infty {(-1)^j r^j j! \over (2k_1,2k_2,2k_1+2k_2+j-1)_j} 
{\ }_2F_1\left( {{k+i(x_1+x_2), k+i(y_1+y_2)}\atop{2k}};r \right)\nn \\
&& \qquad\times p_j(x_1;k_1,k_2-i(x_1+x_2),k_1,k_2+i(x_1+x_2))\nn\\
&& \qquad\times p_j(y_1;k_1,k_2-i(y_1+y_2),k_1,k_2+i(y_1+y_2)).
\label{hahn}
\eea

One can interpret (\ref{hahn}) as a bilinear sum formula for
continuous Hahn polynomials. Using the relabelling $(x,y)=(x_1,y_1)$,
$a=k_1$, $b=k_2-i(x_1+x_2)$, $c=k_2-i(x_1+x_2)$, $b'=k_2-i(y_1+y_2)$,
and $c'=k_2+i(y_1+y_2)$, we have the following
\begin{theo}
Let $a>0$, $\Re(b,c,b',d')>0$ with $b+d=b'+d'$, $\bar d =b$ and $\bar
d'= b'$. Then the continuous Hahn polynomials satisfy
a bilinear sum formula given by~:
\beas
&&\sum_{j=0}^\infty h_j\, p_j(x;a,b,a,d)\, p_j(y,a,b',a,d')\, r^j = \\
&& {}_2F_1\left( {a+ix,a+iy\atop 2a};r \right) 
{\ }_2F_1\left( {d-ix,d'-iy\atop b+d};r\right),
\eeas
where
$$
h_j = {(-1)^j j!\over (2a,b+d,2a+b+d+j-1)_j}
{\ }_2F_1\left( {a+d+j,a+d'+j\atop 2a+b+d+2j};r \right),
$$
and $|r|<1$, $x,y\in\Real$.
\label{theochahn}
\end{theo}

\begin{rema} {\rm
The bilinear generating function of Theorem~\ref{theochahn}
is an extension of
\bea
&&\sum_{j=0}^\infty
{{ (-1)^j\, j!\,(\al+\be+2j+1)\Gamma(\al+\be+j+1)} \over
{\Gamma(\al+j+1) \Gamma (\be+j+1)}}
P_j^{(\al,\be)}(x) P_j^{(\al,\be)}(y)\, J_{\al+\be+2j+1}(z) \nn\\
&& = 2^{\al+\be-1}
\bigl( (1-x)(1-y)\bigr)^{-\al/2}
\bigl( (1+x)(1+y)\bigr)^{-\be/2}\nn\\
&&\qquad \times \, z\,
J_\al \left( {z\over 2} \sqrt{ (1-x)(1-y)}\right)
J_\be \left( {z\over 2} \sqrt{ (1+x)(1+y)}\right).
\label{bilgfjacobi}
\eea
This bilinear generating function for the Jacobi polynomials
$P_n^{(\al,\be)}$ can be obtained from Theorem~\ref{theochahn}
using a standard limit of
the continuous Hahn polynomials to Jacobi polynomials,
cf. \cite[\S 2.8]{KoekS}, and in this rescaling the
${}_2F_1$'s go over
into ${}_0F_1$'s, which can be rewritten in terms of
Bessel functions $J_\nu$ of the first kind \cite{Wats}.
This formula is due to Bateman (1904), see
\cite[\S 11.6]{Wats}.
It should be noted that Bateman's formula
\beq
\sum_{k=0}^n c_{k,n} P_k^{(\al,\be)}(x)P_k^{(\al,\be)}(y)
=(x+y)^n P_n^{(\al,\be)}\Bigl({{1+xy}\over{x+y}}\Bigr),
\label{bateman}
\eeq
where the $c_{k,n}$ follow by specialising $y=1$, can be
obtained from (\ref{bilgfjacobi}) by writing both sides of 
(\ref{bilgfjacobi})
as power series in $z$ after division by $z^{\al+\be+1}$
using \cite[4.3(12)]{HTF1} for the right hand side
and equating coefficients.

In the same way we can obtain from Theorem~\ref{theochahn} an extension
of (\ref{bateman}). }
\end{rema}

\begin{coro}
Assume $b+d=b'+d'$, $\bar d =b$
and $\bar d'= b'$, then
\beas
&&\sum_{j=0}^k {(-k)_j (2a+b+d)_{2j} j! \over
(2a,b+d,2a+b+d+j-1,a+d,a+d',2a+b+d+k)_j}\,\\
&&\qquad\times p_j(x;a,b,a,d)\,p_j(y,a,b',a,d') \\
&& = {(d-ix,d'-iy,2a+b+d)_k \over (a+d,a+d',b+d)_k}
{\;}_4F_3\left( {{-k,1-k-b-d,a+ix,a+iy}\atop{2a,1-k-d+ix,1-k-d'+iy}};
1\right).
\eeas
\label{corohahn}
\end{coro}

Note that for $b=d'$, and thus $d=b'$ and $\bar b=d$,
the ${}_4F_3$-series
is balanced. In this case we can use Whipple's transformation,
see e.g. \cite[(2.10.5)]{GR}, to see that the ${}_4F_3$-series
in Corollary~\ref{corohahn} equals
$$
{{(a+d)_k(a+b+i(x-y))_k}\over{(d-ix)_k(b-iy)_k}}\,
{}_4F_3\left( {{-k,2a+b+d+k-1,a+ix,a-iy}\atop{2a, a+d,
a+b+i(x-y)}};1\right),
$$
which is a Wilson polynomial for $x=y$.

\noindent {\em Proof.} Write the ${}_2F_1$ in $h_j$ in the left hand side
of Theorem~\ref{theochahn} as
$$
\sum_{k=j}^\infty {{(a+d+j)_{k-j}(a+d'+j)_{k-j}}\over
{(k-j)!\, (2a+b+d+2j)_{k-j}}} r^{k-j}
$$
to see that the coefficient of $r^k$ in the left hand side equals
\beas
&&\sum_{j=0}^k {{(a+d+j)_{k-j}(a+d'+j)_{k-j}}\over
{(k-j)!\, (2a+b+d+2j)_{k-j}}} {{(-1)^j j!}\over {(2a,b+d,2a+b+d+j-1)_j}}\\
&&\qquad\qquad\times p_j(x;a,b,a,d)\, p_j(y,a,b',a,d').
\eeas
The coefficient of $r^k$ in the right hand side of 
Theorem~\ref{theochahn} follows from \cite[4.3(14)]{HTF1}.
Equating and rewriting gives the result. \qquad\Box

Formula (\ref{hahn}) can be extended to more general values.
Relabelling the parameters of the first (resp.\ second) ${}_2F_1$ in
(\ref{hahn}) by $a,b,c$ (resp.\ $a',b',c'$), and expressing the
continuous Hahn polynomials in terms of a ${}_3F_2$, we obtain~:

\begin{theo}
For $a,b,c,a',b',c'\in\C$ with $c$ and $c'$ no negative integers, and
$|z|<1$, the following multiplication formula holds~:
$$
{}_2F_1\left( {a,b\atop c};z \right) {\ }_2F_1\left( {a',b'\atop c'};z
\right) = \sum_{j=0}^\infty C_j \, z^j {\ }_2F_1\left( {a+a'+j,b+b'+j
\atop c+c'+2j};z \right),
$$
where
\beas
C_j &=& {(c,a+a',b+b')_j\over j!(c',c+c'+j-1)_j}
{\ }_3F_2\left( {-j,a,c+c'+j-1\atop a+a',c};1 \right)\\
&&\times {\ }_3F_2\left( {-j,b,c+c'+j-1\atop b+b',c};1 \right).
\eeas
\label{theoMP}
\end{theo}

The theorem holds because as a power series in $z$ the coefficient of
$z^k$ in both sides of the equation is the same (this follows from
(\ref{hahn})). Moreover, for $|z|<1$ and $c$ and $c'$ no negative
integers the series are well defined and absolutely
convergent. Note that for
$a+a'$ (or $b+b'$) a negative integer, $C_j$ has to be interpreted with
care and in fact it is not zero for $j>-a-a'$
unless both $a$ and $a'$ are negative integers. In
that case $C_j=0$ for $j>-a-a'$ and we obtain a polynomial
identity. The polynomial identity is equivalent to the dual
statement of the convolution identity 
\cite[Thm.~3.4, Rem.~3.5(ii)]{KV} as is easily
seen by specialising the parameters such that each ${}_2F_1$
in the left hand side of Theorem~\ref{theoMP} corresponds to a
Meixner-Pollaczek polynomial~(\ref{defMP}). Observe that this
result~(\ref{ecoupledMP}) is used in the derivation of 
Theorem~\ref{theoMP}. 

A special case of Theorem~\ref{theoMP} occurs when $(a',b',c')=(a,b,c)$. 
Then the ${}_3F_2$'s in
$C_j$ can be summed by Watson's theorem~\cite[4.4(7)]{HTF1}
and we obtain the Burchnall-Chaundy (1948) formula for the square of
a ${}_2F_1$-series, cf.~\cite[2.5.2(7)]{HTF1}.

An interesting limit of Theorem~\ref{theoMP} is found by putting
$z=x/b$ and $b'=yb/x$, and taking the limit for $b\rightarrow \infty$~:
\begin{coro}
\beas
&&{}_1F_1\left( {a\atop c};x \right) {\ }_1F_1\left( {a'\atop c'};y
\right) = \sum_{j=0}^\infty D_j \, (x+y)^j \\
&&\times {\ }_2F_1\left( {-j,c+c'+j-1
\atop c};{x\over x+y} \right) {\ }_1F_1\left( {a+a'+j\atop c+c'+2j};x+y
\right) , 
\eeas
where
$$
D_j = {(c,a+a')_j\over j!(c',c+c'+j-1)_j}
{\ }_3F_2\left( {-j,a,c+c'+j-1\atop a+a',c};1 \right).
$$
\label{coroMP2}
\end{coro}
Again for $a$ and $a'$ negative integers we obtain the dual
statement of the general convolution identity \cite[(1.1)]{V},
\cite[Cor.~3.6(i)]{KV} for Laguerre polynomials.

As mentioned in the beginning of this section, the other case to study
is the realisation in terms of Meixner polynomials. However, the
analysis is completely similar to that of Meixner-Pollaczek
polynomials (the continuous Hahn polynomials being replaced by Hahn
polynomials), and the results are essentially the same. We mention one
formula here, namely a bilinear generating function for Hahn
polynomials. The definition of the Hahn polynomial reads~\cite{KoekS}
\beq
Q_n(x;\al,\be,N)={}_3F_2\left({-n,n+\al+\be+1,-x \atop \al+1, -N};
1\right), \qquad n=0,1,2,\ldots,N.
\eeq
They satisfy a discrete orthogonality relation, with support
$\{0,1,\ldots ,N\}$. We obtain~:
\begin{theo}
Let $M$ and $N$ be positive integers, $\al,\be,z\in\C$, $x\in\{
0,1,\ldots ,M\}$, and $y\in\{0,1,\ldots, N\}$. Then
\beas
&&\sum_{j=0}^{\min(M,N)} Q_j(x;\al,\be,M) Q_j(y;\al,\be,N) {(\al+1, -M,
-N)_j \over (\be+1, \al+\be+j+1)_j} \\
&&\qquad \times {}_2F_1\left({j-M,j-N\atop \al+\be+2j+2};z\right) z^j
\\ 
&&= {}_2F_1\left({-x,-y\atop \al+1};z\right)
{}_2F_1\left({x-M,y-N\atop \be+1};z\right).
\eeas
\end{theo}

It is easy to see how to obtain this result from Theorem~\ref{theoMP}. 

\section{Polynomials related to $U_q(\mbox{\eultit su}(1,1))$}
\setcounter{equation}{0}

Assume that $0<q<1$.
Let $U_q({\Eu s \Eu l}(2,\C))$ be the complex unital
associative algebra generated by $A$, $B$, $C$, $D$ subject to the
relations
\beq
AD=1=DA, \quad AB=q^{1/2}BA,\quad AC=q^{-1/2}CA,\quad
BC-CB = {{A^2-D^2}\over{q^{1/2}-q^{-1/2}}}.
\label{defU}
\eeq
It is a Hopf algebra with comultiplication
\beq
\begin{array}{l}
\De(A)=A\otimes A,\quad \De(B)=A\otimes B+B\otimes D,\\
\De(C) = A\otimes C+C\otimes D, \quad \De(D)=D\otimes D
\end{array}
\label{defcomult}
\eeq
on the level of generators and extended as an algebra
homomorphism. The $\ast$-structure corresponding to $\U$ is
$$
A^\ast=A,\quad B^\ast=-C,\quad C^\ast=-B, \quad D^\ast=D.
$$

The positive discrete series representations $\pi_k$ of $\U$
are unitary representations labelled by $k>0$. They act in $\Hi$
and the action of the generators is given by
\beq
\begin{array}{ll}
\pi_k(A)\, e^k_n &=q^{(k+n)/2}\, e^k_n, \\[1mm]
\pi_k(C)\, e^k_n &= q^{1/4-(k+n)/2}
{{\sqrt{(1-q^{n})(1-q^{2k+n-1})}}\over{q^{1/2}-q^{-1/2}}}\,
e^k_{n-1},\\[1mm] 
\pi_k(B) \, e^k_n &=
q^{-1/4-(k+n)/2}
{{\sqrt{(1-q^{n+1})(1-q^{2k+n})}}\over{q^{-1/2}-q^{1/2}}}\, e^k_{n+1}.
\end{array}
\label{actionsuq11}
\eeq

Recall that the tensor product of two representations
is defined by use of the comultiplication.
The tensor product of two positive discrete series representations
decomposes as for the Lie algebra $\su(1,1)$, see
eq.~(\ref{tensorproddec}). The Clebsch-Gordan coefficients for $\U$
appearing in
\beq
e^{(k_1k_2)k}_n = \sum_{n_1,n_2} C^{k_1,k_2,k}_{n_1,n_2,n}(q)\,
e^{k_1}_{n_1}\otimes e^{k_2}_{n_2}
\label{defqCGC}
\eeq
are non-zero only if $n_1+n_2=n+j$,
$k=k_1+k_2+j$ for $j,n_1,n_2,n\in\Zp$, and are normalised
by $\langle e^{(k_1k_2)k}_0, e^{k_1}_0\otimes e^{k_2}_j\rangle >0$. In
\cite[Lemma 4.4]{KV} they have been computed in terms of $q$-Hahn
polynomials. 

The above results can be found in Burban and Klimyk \cite{BurbK}
and Kalnins, Manocha and Miller \cite{KalnMM}. See
Chari and Pressley \cite{CharP} for general information on
quantised universal enveloping algebras.

Next, we consider the polynomials related to these representations. Our
notation for $q$-special functions and $q$-shifted factorials is as 
in~\cite{GR}.
The Askey-Wilson polynomials are, see Askey and Wilson~\cite{AW} or
e.g.~\cite[\S 7.5]{GR},
\beq
p_n(x;a,b,c,d|q) = a^{-n} (ab,ac,ad;q)_n\, {\ }_4\vp_3
\left( {{q^{-n},abcdq^{n-1},ae^{i\th},ae^{-i\th}}\atop
{ab,\ ac,\ ad}}; q,q\right),
\label{a-w}
\eeq
with $x=\cos\th$. 
The Al-Salam and Chihara polynomials are obtained by taking
$c=d=0$ in the Askey-Wilson polynomials;
\beq
R_n(x;a,b|q) = p_n(x;a,b,0,0|q) = a^{-n} (ab;q)_n\,
{\ }_3\vp_2 \left( {{q^{-n},ae^{i\th},ae^{-i\th}}\atop{ab,\ 0}};
q,q\right).
\label{a-c}
\eeq
The corresponding normalised polynomials are denoted by
\beq
r_n(x;a,b|q) = R_n(x;a,b|q)/\sqrt{(q,ab;q)_n}.
\eeq
For $a$ and $b$ real, or complex conjugates,
with $\max(|a|,|b|)<1$ the orthogonality measure
$d\mu(x;a,b|q)$ is absolutely continuous on $[-1,1]$ and given by
\beq
{1\over 2\pi} \int_{-1}^{1} (q,ab;q)_\infty {w(x)\over\sqrt{1-x^2}}
r_m(x;a,b|q) r_n(x;a,b|q) dx = \de_{m,n},
\eeq
where
$$
w(x) = {h(x,1)h(x,-1)h(x,q^{1/2})h(x,-q^{1/2}) \over h(x,a)h(x,b)},
$$
with
$$
h(x,\al)=(\al e^{i\th}, \al e^{-i\th};q)_\infty, \qquad x=\cos\th .
$$
If $a>1$, $|b|<1$ and $|ab|<1$, the orthogonality measure contains a
discrete part. Here we shall assume that the parameters are
such that we have an absolute continuous measure; in the final result
this condition can be weakened.

Let $k>0$ and assume that $s\in\Real$ satisfies $q^{k}<|s|<q^{-k}$, or
that $s$ is complex with $|s|=1$. 
There is a unitary operator $\La : \Hi \rightarrow
L^2(\Real,d\mu(x;q^ks,q^ks^{-1}|q)) $ mapping each $e^k_n$ onto
$r_n(x;q^ks,q^ks^{-1}|q)$, and for elements $X\in\U$ their realisation in 
$L^2(\Real,d\mu(x;q^ks,q^ks^{-1}|q))$ will be denoted by $\rh_{k,s}(X)$.

\begin{prop} 
Let $t\in\C,\quad |t|<1$, $x,y\in[-1,+1]$, $|s|,|\si|\in(q^k,q^{-k})$, and
\beq
v^k_t(x,s;y,\si) = \sum_{n=0}^\infty r_n(x;q^ks,q^ks^{-1}|q)
 r_n(y;q^k\si,q^k\si^{-1}|q) t^n. 
\label{vAC}
\eeq
Then\\
(i) $ v^k_t(x,s;y,\si) \in L^2(\Real,d\mu(x;q^ks,q^ks^{-1}|q))$.\\
(ii) Explicitly, with $x=\cos\th$ and $y=\cos\phi$
\bea
&&v^{k}_t(x,s;y,\si)= 
{(q^kte^{-i\phi}s, q^kte^{-i\phi}s^{-1}, 
q^kte^{-i\th}\si,q^kte^{-i\th}\si^{-1} ;q)_\infty \over
(te^{i(\th-\phi)}, te^{-i(\th-\phi)},
te^{-i(\th+\phi)},q^{2k}te^{-i(\th+\phi)} 
;q)_\infty} \nn\\
&&\times{}_8W_7(q^{2k-1}te^{-i(\th+\phi)}; q^k e^{-i\th}s,
q^k e^{-i\th}s^{-1},q^k e^{-i\phi} \si,q^k e^{-i\phi} \si^{-1},
te^{-i(\th+\phi)} ;q, te^{i(\th+\phi)}). \nn\\
&& \ \label{vexplicitAC} 
\eea
Herein, $_{r+1}W_r$ is a very-well-poised $_{r+1}\vp_r$ series, see
\cite[(2.1.11)]{GR}. \\
(iii) Let 
$$
X_t = q^{1/4}tB-q^{-1/4}t^{-1}C + {{\si^{-1}+\si} \over
{q^{-1/2}-q^{1/2}}}(A-D). 
$$ 
Then $v^k_t(x,s;y,\si)$ is an
eigenfunction of $\rh_{k,s}(X_t A)$ for the eigenvalue $(2y-\si-\si^{-1})
/ (q^{-1/2}-q^{1/2})$.
\label{propAC}
\end{prop}

\noindent {\em Proof.} 
Statement (i) follows from the asymptotic behaviour of the Al-Salam
and Chihara polynomials~\cite[\S 3.1]{AI}, see also~\cite{IW,Wilson}.
Expression (\ref{vAC}) is the nonsymmetric
Poisson kernel for the Al-Salam and
Chihara polynomials, and has been deduced by Askey, Rahman and 
Suslov~\cite{ARS}. Using~\cite[(14.5)]{ARS} gives
\bea
&&v^{k}_t(x,s;y,\si)= 
{(t^2,q^k e^{-i\th}s^{-1}, q^kte^{i\th}\si, 
q^kt\si^{-1} e^{i\th},q^ktse^{i\ph},q^ktse^{-i\ph} ;q)_\infty \over
(q^{2k},q^kst^2e^{i\th},te^{i(\th+\phi)}, te^{i(\th-\phi)},
te^{-i(\th+\phi)},te^{-i(\th-\phi)} ;q)_\infty} \nn\\
&&\times{}_8W_7(q^{k-1}st^2e^{i\th}; te^{i(\th+\phi)} ,te^{i(\th-\phi)},
q^kse^{i\th},st\si^{-1},st\si;q,q^ks^{-1}e^{-i\th}). \nn\\
&& \ \label{tmpAC} 
\eea
Using (III.24) of \cite{GR} on this equation yields
(\ref{vexplicitAC}). The rest of the statements are proved in a
similar way as Proposition~\ref{propMP}, using in this
case the three-term recurrence relation for Al-Salam and Chihara
polynomials. \qquad\Box

In the tensor product $\pi_{k_1}\otimes \pi_{k_2}$,
there is a unitary operator $\Upsilon : \Hi\otimes\Hi \rightarrow L^2(
\Real^2, d\mu(x_1;q^{k_1}s_1,q^{k_1}s_1^{-1}|q) 
d\mu(x_2;q^{k_2}s_2,q^{k_2}s_2^{-1}|q) )$ mapping
$e^{k_1}_{n_1} \otimes e^{k_2}_{n_2}$ onto
$r_{n_1}(x_1;q^{k_1}s_1,q^{k_1}s_1^{-1}|q) 
r_{n_2}(x_2;q^{k_2}s_2,q^{k_2}s_2^{-1}|q) $.
Here we shall only consider the case $s_2=s$ and $s_1 =e^{i\th_2}$
($x_2= \cos\th_2$) and denote the corresponding measure by
$$
d\mu(x_1,x_2)=d\mu(x_1;q^{k_1}e^{i\th_2},q^{k_1}e^{-i\th_2}|q) 
d\mu(x_2;q^{k_2}s,q^{k_2}s^{-1}|q) .
$$
Then we can simplify the basis
functions in the tensor product decomposition. The realisation of
$\Delta(X)$ with $X\in \U$ in $L^2(
\Real^2, d\mu(x_1,x_2))$ is denoted by 
$\si_{k_1,k_2,s}(X)$.
For $k=k_1+k_2+j$, with $j\in\Zp$, and $x_i=\cos\th_i$, 
the above-mentioned simplification reads~\cite[(4.11)]{KV}
\bea
e^{(k_1k_2)k}_n(x_1,x_2;s) &=& \sum_{n_1,n_2} C^{k_1,k_2,k}_{n_1,n_2,n}(q)\
 r_{n_1}(x_1;q^{k_1}e^{i\th_2},q^{k_1}e^{-i\th_2}|q) 
 r_{n_2}(x_2;q^{k_2}s,q^{k_2}s^{-1}|q)  \nn\\
&=& r_{n}(x_1;q^ks,q^ks^{-1}|q)\, S^{(k_1,k_2)}_j(x_1,x_2;s),
\label{ecoupledAC}
\eea
where in this case $S_j$ is given in terms of an Askey-Wilson
polynomial~: 
\beq
S^{(k_1,k_2)}_j(x_1,x_2;s)={p_j(x_2;q^{k_1}e^{i\th_1},q^{k_1}e^{-i\th_1},
q^{k_2}s, q^{k_2}s^{-1}|q) \over \sqrt{(q,q^{2k_1},q^{2k_2},q^{2k_1
+2k_2+j-1}; q)_j}}.
\label{SexplicitAC}
\eeq
Now we have the following
\begin{prop}
Let $t\in\C,\quad |t|<1$, $x_1,x_2,y\in[-1,+1]$, 
$|s|,|\si|\in(q^k,q^{-k})$, and
\beq
w^{(k_1k_2)k}_t(x_1,x_2,s;y,\si) = \sum_{n=0}^\infty
r_n(y;q^k\si,q^k\si^{-1}|q)\,
 e^{(k_1k_2)k}_n(x_1,x_2;s)\, t^n
\label{wAC}
\eeq
for $k=k_1+k_2+j$ with $j\in\Zp$.
Let again $X_t$ be as in Proposition~\ref{propAC}(iii), then\\
(i) $w^{(k_1k_2)k}_t(x_1,x_2,s;y,\si) \in L^2(\Real^2,d\mu(x_1,x_2))$
is an eigenfunction of
$\si_{k_1,k_2,s} (X_t A)$ for the eigenvalue  $(2y-\si-\si^{-1})
/ (q^{-1/2}-q^{1/2})$.\\
(ii) Explicitly, one has 
\beq
w^{(k_1k_2)k}_t(x_1,x_2,s;y,\si) = S_j^{(k_1,k_2)}(x_1,x_2;s)\,
v^k_t(x_1,s;y,\si). 
\label{wexplicitAC}
\eeq
(iii) For $-1\leq y_1,y_2 \leq 1$, 
\beq
v^{k_1}_t(x_1,e^{i\th_2};y_1,e^{i\ph_2})v^{k_2}_t(x_2,s;y_2,\si) =
\sum_{j=0}^\infty t^j S_j^{(k_1,k_2)}(y_1,y_2;\si)\,
w^{(k_1k_2)k}_t(x_1,x_2,s;y_1,\si),
\label{CGwAC}
\eeq
where $x_i=\cos\th_i$ and $y_i=\cos\ph_i$.
\label{propwAC}
\end{prop}

In this case, the basic identity follows from (\ref{wexplicitAC}) and
(\ref{CGwAC})~: 
\bea
&&v^{k_1}_t(x_1,e^{i\th_2};y_1,e^{i\ph_2})v^{k_2}_t(x_2,s;y_2,\si) =
 \nn\\
&&\sum_{j=0}^\infty t^j v^{k_1+k_2+j}_t(x_1,s;y_1,\si)\, 
S_j^{(k_1,k_2)}(x_1,x_2;s)\, S_j^{(k_1,k_2)}(y_1,y_2;\si) .
\label{SpoissonAC}
\eea
Using the explicit expressions (\ref{vexplicitAC}) and
(\ref{SexplicitAC}), relabelling $(x_2,y_2)$ by $(x,y)$ and putting
$$
\begin{array}{ll}
(a,b,c,d)&=(q^{k_1}e^{i\th_1},q^{k_1}e^{-i\th_1},q^{k_2}s,q^{k_2}s^{-1}),\\
(a',b',c',d')&=(q^{k_1}e^{i\ph_1},q^{k_1}e^{-i\ph_1},q^{k_2}\si,
q^{k_2}\si^{-1}) ,
\end{array}
$$
we have the following
\begin{theo}
For $|t|<1$, $(a,b,c,d)$ with
$\max(|a|,|b|,|c|,|d|)<1$, $(a',b',c',d')$ with
$\max(|a'|,|b'|,|c'|,|d'|)<1$, $|a't/b|<1$, and 
\beq
ab=a'b',\qquad cd=c'd',
\label{condition}
\eeq
we have, with $x=\cos\th$ and $y=\cos\ph$,
\bea
&& \sum_{j=0}^\infty H_j\, p_j(x;a,b,c,d|q) p_j(y;a',b',c',d'|q) t^j =
\nn\\ 
&&{(bte^{i\ph},bte^{-i\ph},cte^{-i\ph},dte^{-i\ph},
 b'te^{i\th},b'te^{-i\th},c'te^{-i\th},d'te^{-i\th};q)_\infty  \over
(bb't,te^{i(\th-\ph)},te^{i(\ph-\th)},te^{-i(\th+\ph)},
cdte^{-i(\th+\ph)};q)_\infty} \nn\\
&& \times {}_8W_7(bb't/q;be^{i\th},be^{-i\th}, b'e^{i\ph},b'e^{-i\ph},
bt/a' ; q, a't/b) \nn\\
&& \times {}_8W_7(cdte^{-i(\th+\ph)}/q; ce^{-i\th}, de^{-i\th},
c'e^{-i\ph}, d'e^{-i\ph}, te^{-i(\th+\ph)}; q, te^{i(\th+\ph)}),
\label{resultAW}
\eea
where
\beas
&&H_j = {(bc'q^jt,b'cq^jt, bd'q^jt,b'dq^jt;q)_\infty \over
(bb'cdq^{2j}t;q)_\infty (q,ab,cd,abcdq^{j-1};q)_j } \\
&&\times {}_8W_7(bb'cdq^{2j-1}t; bcq^j, bdq^j, b'c'q^j, b'd'q^j, bt/a';
q, a't/b ).
\eeas
\label{theoAW}
\end{theo}

The proof follows directly from (\ref{SpoissonAC}). The result
(\ref{resultAW}) holds for arbitrary parameters satisfying
(\ref{condition}) as long as all the series in (\ref{resultAW}) converge
absolutely. The absolute convergence follows from the fact that for 
$j\rightarrow\infty$ the coefficients $H_j$ go to 
$(t^2;q)_\infty / (q,ab,cd,a't/b;q)_\infty$, together with the
asymptotic behaviour of Askey-Wilson polynomials~\cite{IW,Wilson}.
The condition $|a't/b|<1$ is not essential as one can use Bailey's
transform \cite[(III.23-24)]{GR} to rewrite the ${}_8W_7$-series.

Rahman \cite[(4.4)]{Rahm2} derives a formula remotely similar
to (\ref{resultAW}) in the sense that the summand also
consists of the product of a ${}_8W_7$-series and the 
product of two Askey-Wilson polynomials. However, the right
hand side in (\ref{resultAW}) is much less complicated than the 
right hand side of \cite[(4.4)]{Rahm2}, which consists of four
infinite sums each involving ${}_{10}W_9$-series. Rahman's
formula is a $q$-analogue of Feldheim's bilinear generating
function for the Jacobi polynomials, whereas (\ref{resultAW})
is a $q$-analogue of Bateman's bilinear generating
function for the Jacobi polynomials, see (\ref{bilgfjacobi}).
It should be interesting to investigate whether (\ref{resultAW}) 
could lead to a $q$-analogue of (\ref{bateman}) or 
Corollary~\ref{corohahn}.

The previous formulas can be specialised by putting $d=d'=0$; then the
Askey-Wilson polynomials become continuous dual $q$-Hahn polynomials
denoted by $p_n(x;a,b,c|q)$. One finds~:
\begin{coro}
For given $|t|<1$, $(a,b,c)$ and $(a',b',c')$ with
$\max(|a|,|b|,|c|)<1$ and $\max(|a'|,|b'|,|c'|)<1$, $|a't/b|<1$,
and
\beq
ab=a'b',
\label{condition2}
\eeq
we have, with $x=\cos\th$ and $y=\cos\ph$,
\bea
&& {\sum_{j=0}^\infty} G_j\, p_j(x;a,b,c|q) p_j(y;a',b',c'|q) t^j =
\nn\\ 
&&{(bte^{i\ph},bte^{-i\ph},cte^{-i\ph},
 b'te^{i\th},b'te^{-i\th},c'te^{-i\th};q)_\infty  \over
(bb't,te^{i(\th-\ph)},te^{i(\ph-\th)},te^{-i(\th+\ph)};q)_\infty} \nn\\
&& \times {}_8W_7(bb't/q;be^{i\th},be^{-i\th}, b'e^{i\ph},b'e^{-i\ph},
bt/a' ; q, a't/b) \nn\\
&& \times {}_3\vp_2\left( {ce^{-i\th}, c'e^{-i\ph}, te^{-i(\th+\ph)}
\atop cte^{-i\ph},c'te^{-i\th}}; q, te^{i(\th+\ph)}\right),
\label{resultqH}
\eea
where
\beq
G_j = {(bc'q^jt,b'cq^jt;q)_\infty \over (q,ab;q)_j } 
{\ }_3\vp_2\left({bcq^j, b'c'q^j, bt/a' \atop bc'q^jt, b'cq^jt};
q, a't/b \right).
\eeq
\label{coroqH}
\end{coro}

By making the further specialisation $c=c'=0$ the identity reduces to
the non-symmetric Poisson kernel for the Al-Salam and Chihara
polynomials~\cite[(14.8)]{ARS}, i.e.\ essentially the formula that we
started from and used in~(\ref{vexplicitAC}).

Another specialisation of Corollary~\ref{coroqH} can be performed as 
follows. In $G_j$, use transformation~\cite[(III.9)]{GR} for the 
${}_3\vp_2$-series; on the right hand side
of (\ref{resultqH}) use 
transformation~\cite[(III.23)]{GR} for the ${}_8W_7$-series with
$(b,c,d)$ corresponding to $(be^{-i\th}, b'e^{-i\phi}, bt/a)$; finally 
multiply both sides of  (\ref{resultqH}) by  $(a't/b;q)_\infty$ and
take $b=b'=0$. The result is a bilinear generating function for
Al-Salam and Chihara polynomials without the extra assumption that
the product of the parameters should be equal.

\begin{coro} 
For $|t|<1$ and
$\max(|a|,|c|)<1$, $\max(|a'|,|c'|)<1$ and with
$x=\cos\th$, $y=\cos\phi$ we have
\beas
&&\sum_{j=0}^\infty {{t^j}\over{(q,a'ct;q)_j}}
\, {}_2\vp_1\left( {{ct/c',acq^j}\atop{a'ctq^j}};q, t{c'\over c}
\right) \ R_j(x;a,c\mid q) R_j(y;a',c'|q) = \\
&&{{(cte^{-i\phi},c'te^{-i\th},ate^{i\phi},a'te^{i\th};q)_\infty}
\over{(te^{i(\th-\phi)},te^{i(\phi-\th)},c't/c,a'ct;q)_\infty}}
{}_3\vp_2\left( {{a'e^{i\phi},ae^{i\th}, te^{i(\th+\phi)}}\atop{
ate^{i\phi}, a'te^{i\th}}};q,te^{-i(\phi+\th)}\right) \\
&&\qquad\qquad\qquad
\times {}_3\vp_2\left( {{ce^{-i\th},c'e^{-i\phi}, te^{-i(\th+\phi)}}
\atop{cte^{-i\phi}, c'te^{-i\th}}};q,te^{i(\th+\phi)}\right).
\eeas
\label{coroAC}
\end{coro}

\begin{rema}
{\rm Take $a=a'=0$ in Corollary~\ref{coroAC}, so that the
${}_2\vp_1$-series and the first ${}_3\vp_2$-series reduce
to ${}_1\vp_0$-series which can be summed by the $q$-binomial
theorem \cite[(II.3)]{GR}. The resulting formula is the
non-symmetric Poisson kernel for the continuous big $q$-Hermite
polynomials, i.e. Askey-Wilson polynomials with three parameters
set to zero. So for $a=a'=0$ Corollary~\ref{coroAC} overlaps with
Askey, Rahman and Suslov \cite[(14.8)]{ARS}. }
\end{rema}

\vskip 3mm
{\sc Vakgroep Wiskunde, Universiteit van Amsterdam, Plantage Muidergracht
24, 1018 TV Amsterdam, The Netherlands}

{\em E-mail address}~: koelink@wins.uva.nl

\vskip 3mm \noindent
{\sc Vakgroep Toegepaste Wiskunde en Informatica, Universiteit Gent,
Krijgslaan 281-S9, B-9000 Gent, Belgium}

{\em E-mail address}~: Joris.VanderJeugt@rug.ac.be

\end{document}